\documentclass[12pt]{article}
\usepackage{amsmath, latexsym, amssymb}
\input xypic

\begin{document}

\title{Weak equivalence classes of complex vector bundles }
\author{\bf  H\^ong-V\^an L\^e \\
Mathematical Institute  of ASCR\\
Zitna 25, Praha 1, CZ-11567\\
and\\
Max-Planck-Institute f\"ur Mathematik\\
Inselstr.22, D-04103 Leipzig}

\date{}
\maketitle
\newcommand{\R}{{\mathbb R}}
\newcommand{\C}{{\mathbb C}}
\newcommand{\F}{{\mathbb F}}
\newcommand{\Z}{{\mathbb Z}}
\newcommand{\N}{{\mathbb N}}
\newcommand{\Q}{{\mathbb Q}}
\newcommand{\Hq}{{\mathbb H}}

\newcommand{\Aa}{{\mathcal A}}
\newcommand{\Bb}{{\mathcal B}}
\newcommand{\Cc}{{\mathcal C}}    
\newcommand{\Dd}{{\mathcal D}}
\newcommand{\Ee}{{\mathcal E}}
\newcommand{\Ff}{{\mathcal F}}
\newcommand{\Gg}{{\mathcal G}}    
\newcommand{\Hh}{{\mathcal H}}
\newcommand{\Kk}{{\mathcal K}}
\newcommand{\Ii}{{\mathcal I}}
\newcommand{\Jj}{{\mathcal J}}
\newcommand{\Ll}{{\mathcal L}}    
\newcommand{\Mm}{{\mathcal M}}    
\newcommand{\Nn}{{\mathcal N}}
\newcommand{\Oo}{{\mathcal O}}
\newcommand{\Pp}{{\mathcal P}}
\newcommand{\Qq}{{\mathcal Q}}
\newcommand{\Rr}{{\mathcal R}}
\newcommand{\Ss}{{\mathcal S}}
\newcommand{\Tt}{{\mathcal T}}
\newcommand{\Uu}{{\mathcal U}}
\newcommand{\Vv}{{\mathcal V}}
\newcommand{\Ww}{{\mathcal W}}
\newcommand{\Xx}{{\mathcal X}}
\newcommand{\Yy}{{\mathcal Y}}
\newcommand{\Zz}{{\mathcal Z}}

\newcommand{\zt}{{\tilde z}}
\newcommand{\xt}{{\tilde x}}
\newcommand{\Ht}{\widetilde{H}}
\newcommand{\ut}{{\tilde u}}
\newcommand{\Mt}{{\widetilde M}}
\newcommand{\Llt}{{\widetilde{\mathcal L}}}
\newcommand{\yt}{{\tilde y}}
\newcommand{\vt}{{\tilde v}}
\newcommand{\Ppt}{{\widetilde{\mathcal P}}}

\newcommand{\Remark}{{\it Remark}}
\newcommand{\Proof}{{\it Proof}}
\newcommand{\ad}{{\rm ad}}
\newcommand{\Om}{{\Omega}}
\newcommand{\om}{{\omega}}
\newcommand{\eps}{{\varepsilon}}
\newcommand{\Di}{{\rm Diff}}
\newcommand{\Pro}[1]{\noindent {\bf Proposition #1}}
\newcommand{\Thm}[1]{\noindent {\bf Theorem #1}}
\newcommand{\Lem}[1]{\noindent {\bf Lemma #1 }}
\newcommand{\An}[1]{\noindent {\bf Anmerkung #1}}
\newcommand{\Kor}[1]{\noindent {\bf Korollar #1}}
\newcommand{\Satz}[1]{\noindent {\bf Satz #1}}

\newcommand{\gl}{{\frak gl}}
\renewcommand{\o}{{\frak o}}
\newcommand{\so}{{\frak so}}
\renewcommand{\u}{{\frak u}}
\newcommand{\su}{{\frak su}}
\newcommand{\ssl}{{\frak sl}}
\newcommand{\ssp}{{\frak sp}}

\newcommand{\Cinf}{C^{\infty}}
\newcommand{\CS}{{\mathcal{CS}}}
\newcommand{\YM}{{\mathcal{YM}}}
\newcommand{\Jreg}{{\mathcal J}_{\rm reg}}
\newcommand{\Hreg}{{\mathcal H}_{\rm reg}}
\newcommand{\SP}{{\rm SP}}
\newcommand{\im}{{\rm im}}

\newcommand{\inner}[2]{\langle #1, #2\rangle}    
\newcommand{\Inner}[2]{#1\cdot#2}
\def\NABLA#1{{\mathop{\nabla\kern-.5ex\lower1ex\hbox{$#1$}}}}
\def\Nabla#1{\nabla\kern-.5ex{}_#1}

\newcommand{\half}{\scriptstyle\frac{1}{2}}
\newcommand{\p}{{\partial}}
\newcommand{\notsub}{\not\subset}
\newcommand{\iI}{{I}}               
\newcommand{\bI}{{\partial I}}      
\newcommand{\LRA}{\Longrightarrow}
\newcommand{\LLA}{\Longleftarrow}
\newcommand{\lra}{\longrightarrow}
\newcommand{\LLR}{\Longleftrightarrow}
\newcommand{\lla}{\longleftarrow}
\newcommand{\INTO}{\hookrightarrow}

\newcommand{\Sy}{\text{ Diff }_{\om}}
\newcommand{\Ex}{\text{Diff }_{ex}}
\newcommand{\jdef}[1]{{\bf #1}}
\newcommand{\QED}{\hfill$\Box$\medskip}

\newcommand{\UuU}{\Upsilon _{\delta}(H_0) \times \Uu _{\delta} (J_0)}
\newcommand{\bm}{\boldmath}
\medskip

\abstract
For any complex vector bundle $E^k$  of rank $k$ over a manifold $M^m$  with  Chern classes
$c_i \in H^{2i}(M^m,\Z)$  and any non-negative integers $l_1, \cdots, l_k$ we show the existence of a positive number
$N(k,m)$ and the existence of a complex vector bundle $\hat E^k$  over $M^m$ whose Chern classes are
$ N(k,m) \cdot l_i\cdot  c_i\in H^{2i} (M^m,\Z)$. We also discuss a  version of this statement for
holomorphic vector bundles  over  projective algebraic manifolds.
\endabstract

\medskip


\section {Introduction.}

The study of complex vector bundles of rank $k$ over a manifold
$M^m$ can be reduced to the study of mappings from $M ^m$ to the
classifying space $Gr_k (\C^\infty) = BU_ k$. Certain equivalence
relations of complex vector bundles lead us  to  study {\bf stable
mappings} of $BU_{k}$ to itself.

We call {\bf a map $g: BU_k \to BU_k$    stable}, if the
restriction of $g$ to any  subspace $Gr_k (\C^N)$ sends $Gr_k
(\C^N)$ to  some Grassmannian $Gr_k (\C^{f (N)})$, moreover
$g^*(c_k) =\lambda \cdot c_k$ for some positive $\lambda$. Here $c_k$
denotes the top Chern class of the universal bundle over $BU_k$.

Two  maps  $f_1, f_2 : M ^m \to  BU_k$  are said to be in one {\bf stable
equivalence class}, if $f_1  = g\circ f_2$ for some stable map $g: BU_k
\to BU_i$. Two complex vector bundles $E_1 ^k$ and $E_2 ^k$ are
said to be in the same {\bf weak equivalence class}, if
the corresponding homotopy classes of
classifying maps contain maps in the
same stable equivalence class. Two complex vector bundles $E_1 ^k$ and
$E_2 ^k$ are called \textbf{ Chern weakly  equivalent}, if  their top Chern classes are differed by a positive constant. Clearly vector bundles in the same weak equivalence class are Chern weakly  equivalent.  Zero sections of 
Chern weakly  equivalent  vector bundles  realize the same homology classes up to a
positive constant.

\medskip

{\bf 1.1. Theorem.} {\it For any complex vector bundle $E^k$  of
rank $k$ over a manifold $M^m$  with the Chern classes $c_i \in
H^{2i}(M^m,\Z)$  and any non-negative integers $l_1, \cdots, l_k$,
$l_k >0$,  there exists a vector bundle $\hat E ^k$ in the same weak equivalence class with
$E^k$, and a positive number $p(k,m)$ such that the Chern
classes $c_i (\hat E ^k)$ are $ N(k,m)\cdot l_i \cdot c_i\in H^{2i}
(M,\Z)$.}
\medskip

As a corollary of Theorem 1.1 and  Thom's theorem [Thom1954,
Theorem II.25] (a detailed proof of this theorem is given in
[Le2005b] and in the Appendix  below) we get

\medskip

{\bf  1.2. Corollary.} [Le2005, Proposition 2.7] {\it  Suppose that $M^m$ is an orientable differentiable manifold.
For any  $c \in  H^{2k} (M ^m, \Z)$ there exists a number $N > 0$ such that there exists a complex vector bundle $E^k$
of rank $k$ over  $M$ whose top Chern class is $N\cdot c$ and all other lower Chern classes are zero.}

\medskip

We can think of Theorem 1.1 together with the Thom theorem  as a
version of the Atiyah-Hirzebruch theorem about isomorphism between the
two rings $K (M^m) \otimes \Q$ and  $H^{even}(M^m, \Q)$ via the
Chern character [A-H1961] which implies that giving an element of
$K(M^m) \otimes \Q$ is the same as giving an element in $H^{even}
(M^m,\Q)$.  Our Theorem 1.1 concerns vector bundles with a given
dimension on $M^m$.  I did not give enough details to the proof of
Proposition 2.7 in [Le2005], so now the proof of Theorem 1.1
should compensate that  deficit.

\medskip

 For $i\le k$  there is a projection map $p$ from $B_{U(i)}$ to $B_{U(k)}$ with fiber
 $U(k)/U(i)$ such that $p^* (c_i) = c_i$. Another consequence of a proof of Theorem 1.1 is

\medskip

{\bf 1.3. Corollary}.  {\it There are maps $g_{k,i}^m : Gr_k
(\C^m) \to Gr_i(\C ^l) \to B_{U(i)}\to  B_{U(k)}$ such that
$g_{k,i}^N (c_j) =p(k,{m\choose k})\cdot \delta^i _j \cdot c_j$ for $1\le i,j \le k$.}

\medskip

In the third section of this note we discuss the problem of extending the
notion of weak equivalence to the category of holomorphic vector
bundles over projective algebraic manifolds.
Two holomorphic vector bundles $E^k$ and $F^k$ over a complex manifold $M^m$ are said to be  {\bf
K\"ahler weakly equivalent},  if there  there are two holomorphic line
bundles $L_1$ and $L_2$ over  $M^m$ such that $E^k\otimes L_1$ and $F^k \otimes L_2$
are  Chern weakly   equivalent.

 It is well-known that the Hodge conjecture is equivalent
to the statement that the Hodge group $H^{p,p}(M,\Q): = H^{2p} (M
,\Q) \cap H^{p,p} (M, \C)$ is generated by the top Chern classes
of holomorphic vector bundle of rank $p$ on a projective algebraic
manifold $ M^n $ (see e.g. [Voisin2002]). A motivation for the
notion of K\"ahler weak equivalence is   Lemma 3.1 below which
states that any holomorphic vector bundle on a projective
algebraic is  K\"ahler weak equivalent to  a holomorphic vector bundle such that the homotopy class of its classifying maps contains a
 holomorphic map.  With this on hand, we  speculate about reduction of the Hodge conjecture
 to the existence of certain holomorphic maps which may be
 obtained by using on the one hand  Zucker  and Saito  results on the existence  of
 normal functions associated to  primitive Hodge cocycles in middle dimensions and on the other hand Siu's technique
 on harmonic maps. 

For the convenience of the reader I include in this note an
appendix which re-exposes the detailed  proof of Thom's theorem in
[Le2005b], which  now has  a simpler form, since this proof is very close to our proof of Theorem 1.1.

\medskip

\section {Proof of  Theorem 1.1.}

{\it Proof of Theorem 1.1.} Denote by $\gamma^k$ the  tautological
bundle over the Grassmannian $Gr_k (\C^N)$ (we  assume that $N
=\infty$ or $N$ is sufficiently large as it shall be specified
later). Since $E^k$ is the pull-back of $\gamma$ via a classifying
map $f : M^m \to Gr_k (\C^N)$, it suffices to   prove Theorem 1.1
for the case $M = Gr_k (\C^N)$, $E^k= \gamma ^k$ and $N$ is
sufficiently large, and after that we use the classifying map $f$ to take back
the  onbtained bundle to $M^m$. 

Let us denote by $K(\Z, n)$ the Eilenberg-McLane  space and by
$\tau^n$ the fundamental class of $K(\Z, n)$. Let $f_k^N : Gr_k
(\C^N) \to K(\Z, 2k)$ be  a classifying map for $c_k (\gamma) \in
H^{2k} (M,\Z)$, i.e.  $(f_k^N)^*(\tau^{2k}) = c_k (\gamma) \in
H^{2i} (M,\Z)$.

\medskip

Let $F_N^n$ be  a map  from $K (\Z, n) \to K(\Z, n)$ such that
$F_N^n (\tau^n) = N \tau^n$. The existence of a map $F_N ^n$ is
ensured  by  the fact that $K(\Z,n)$ is  the classifying space for
$(H^*n, \Z)$. Clearly $F_N^n$ is defined uniquely up to homotopy.

\medskip

{\bf 2.1. Lemma.} [Thom1954, Lemma II.22] {\it  For any finite abelian group $G$ of order $N$
the endomorphism $(F_N^n) ^* : H^* (K(\Z, n), G)$ is trivial.}

\medskip

(Lemma 2.1  follows directly from Cartan's result  which states
that the algebra $H^* (K(\Z, n), \Z_p)$ is generated by iteration
of the Steenrod  squares of  $\tau ^n$).

\medskip

Denote by $Y^q$ the q-skeleton of a CW-complex $Y$. Clearly $\pi_k (K^q (\Z, n)) =
\pi_k (K(\Z,n))$ for any $k\le q$.

\medskip

{\bf 2.2. Proposition.}   {\it Suppose that $Y$ is a  simplicial
space whose q-skeleton $Y^q$ is compact for each $q$. Let the free
component of $\pi_k (Y)$ is isomorphic to $\Z$ with a generator
$t$ and let $Q$ be an integer such that $Q\ge k$.  If for all
$Q\ge q \ge k$ the group  $H^{q+1} (K(\Z, k), \pi_q (Y))$ is
finite, then there exists a map $G_Q : K^Q(\Z, k) \to Y$ such that
$(G_Q)_* ( \pi_k ( K^Q(\Z, k))) = <N(Q,k) t>_{\otimes \Z}\subset
\pi_k (Y)$.}

\medskip

Proposition 2.2 is  a reformulation of Lemma II.24 in [Thom1954],
where Thom did not  explicitly introduce the parameter $Q$.  We
quickly recall his argument, adapted  to this new reformulation.
We prove Proposition 2.2 by induction on the dimension $Q\ge k$.
Clearly Proposition 2.2 for $Q= k$ is trivial, since $K^k (\Z, k)
= S^k$.

Suppose that  we have constructed a map $G_Q$ for $Q\ge k$.

Now we put
$$G_Q^1= F_N ^k \circ G_Q,$$
where $F_N^k$ is the map in Lemma 2.1. By theorem of simplicial approximation we can assume  that $F_N^k$ sends
$K^q (\Z, k)$ to $K^q (\Z, k)$ for each $Q \ge q \ge k$.

Since the obstruction to an extension of $G_Q^1$ to $K^{Q+1} (\Z, k)$
lies in the group $F_N ^* ( H ^{q+1}(K(\Z, k), \pi_q (Y))$ which is trivial by Lemma 2.1, we  shall put $G_{Q+1}$ as  an extension of $G_q^1$ to $K^{Q+1}(\Z, k)$. This completes the induction step for the proof of Proposition 2.2.

\medskip

{\bf 2.3. Lemma.} {\it Suppose that $N \ge 2l+1$. Then   the space $Gr_l (\C^N)$ satisfies the condition for $Y$ with
$Q = N$  and  for all $k = 2r$, if $1\le r\le l$, in Proposition 2.2}.

\medskip

{\it Proof.} To prove Lemma 2.3 it suffices to verify the
following three identities
$$\pi _{2r}  ( Gr_l (\C^N))\otimes \Q =  \Q,\text { for all } 1\le r\le l \leqno (2.3.1)$$
$$\pi _q (Gr_l (\C ^N)) \otimes \Q = 0,  \: \text{ for all other } q\le N \leqno (2.3.2)$$
$$ H^{q+1} ( K (\Z, 2l), \pi_q (Gr_l (\C^N)) \otimes  \Q = 0,\, \forall   q. \leqno (2.3.3)$$

To prove (2.3.1) we consider the following exact sequence
$$\pi_q (U_l \times U_{N-l}) \to \pi_q ( U_N ) \to \pi_q(Gr _l (\C^N)) \to \pi_{q-1} (U_l \times U_{N-l})\leqno (2.3.4)$$
which also remains exact after tensoring with $\Q$. To save the notation  we shall consider this exact sequence as of
that of rational homotopy groups.

\medskip

To prove (2.3.2) we have to consider several cases for $q$.  For
$2\le q = 2r\le 2l$ the exact sequence (2.3.4) implies the
equality (2.3.1), since $\pi_{2r} (U_N)\otimes \Q = 0$,
$\pi_{2r-1} ( U_m)\times \Q = \Q$ for $r\le m$, and taking into
account  that the kernel of the map
$$i: \Q \oplus \Q = \pi_{2r-1} (U_l \times U_{N-l})\otimes \Q \to \pi_{2r-1} (U_N)\times \Q  = \Q$$
 is equal to $\Q$.

\medskip

Now let  us consider the exact sequence (2.3.4)  for $2l +1 \le q \le N$. We know [Spanier1966, 9.7] that
$\pi_q (U_l \otimes U_{N-l}) \otimes \Q = \pi_q ( U_{N-l}) \otimes Q $ and taking into account the fact that  the map
$$ i : \Q = \pi _q (U_l\times U_{N-l} )\otimes \Q \to \pi_q (U_N)\otimes \Q$$
is isomorphism. Taking into account the fact that $\pi_q
(U_{N-l})\otimes Q$ vanishes if $q$ is even, we get
$$ \pi_q (Gr_l (\C^N)) = \ker ( \pi_{q-1} (U_l \times U_N) \to \pi_{q-1} (U_N))= 0$$
 which implies (2.3.2) for $2l+1\le q \le N$.

 \medskip

Finally to verify (2.3.2) for $q$ odd and less than  $2l$ we
notice that the map $\pi_q (U_l\times U_{N-l} ) \to \pi_q (U_N)$
is surjective, hence $\pi_q(Gr_l (\C^N)) = \pi_{q-1} (U_l \times
U_{N-l}) = 0$.

\medskip

The last statement (2.3.3) follows from (2.3.1) for $q = 2l$, and
it follows from (2.3.2)  for all other and taking into account the
fact that $H^* (K(\Z, 2l), \Q) = \Q [x], \, \dim x = 2l$. The last
fact is obtained by Serre and Cartan (see e.g.[F-F, 3.25] for an
exposition. In fact this computation of $H^*(K(\Z, 2l), \Q)$ can
be easily obtained by using induction method and by using the
cohomology spectral sequence associated with the fibration $K(\Z,
n-1) \cong \Om K(\Z, n) \to K(\Z, n)$, whose fiber is
contractible.)\QED

\medskip

{\it Continuation of the proof of Theorem 1.1.}

For $N \ge 2k+1$  and for all $2\le i \le k$ Proposition 2.2 and Lemma 2.3 give us a map
$$G^N_{k,i} : K^{N} (\Z, 2i) \to  Gr_k (\C^N)$$
such that $(G^N_{k,i})_* (w_i) = \alpha (N,k,i) t_i$, where $w_i$
is a generator of $\pi_{2k} (K^N(\Z, 2i)) = \Z$ and $t_i$ is a
generator of $\pi_i (Gr_k (\C^N))\otimes \Q$. Since $H^* (Gr_k
(\C^N),\Z)$ is generated by $c_i (\gamma), i = \overline{1, k},$
[Borel1953], we have
$$ < c_i, t_i>  = A_i \not= 0\leqno(2.4)$$
because $t_i$ is the generator of the free part of $\pi_{2i} (Gr_k (\C^N))$.
(To show that $A_i\not = 0$ we consider the exact sequence (2.3.4). We see easily that the image of $\rho(t_i)$
via embedding $G_k (\C^N)\to G_k (\C^\infty)$ is also a generator of $\pi_i ( Gr_k (\C^\infty))$. Applying the
$\Cc$-version of the Whitehead theorem [Serre1953, Theorem III.3] to $BU_k$ and  the product  $K(\Z, 2)\times \dots \times K(\Z, 2k)$
we notice that
$$ < c_i, i(t_i)>  = A_i \not = 0$$
which  pull back to $G_k (\C^N)$   must also hold.)
Thus
$$ (G^N_{k,i})^* (c_i) = \alpha(N,k,i) \cdot A_i \cdot \tau^{2i}. \leqno (2.5)$$

We can assume that $A_i$ is positive by choosing appropriate orientation of the generator $t_i$.

\medskip

{\it Completion of the proof of  Theorem 1.1.}

Denote by $\lambda _{k,i}^N$ the classifying map from
$G_k(\C^N)$ to $K(\Z, 2i)$ for $c_i\in H^{2i} (G_l (\C^N),\Z)$, i.e.
$$ \lambda_{k,i}^N (\tau ^{2i}) = c_i.$$
We can assume that $\lambda _{k,i}^N ( G_k (\C^N) ) \subset  K^{2k (N-k)} (\Z, 2i)$.
For each $i$ denote by $s (2i)$ the smallest positive number such that  for any $j\le i-1$ and any
$c\in H^{2j} (K (\Z, 2i), \Z)$  we have $s(2i)\cdot  z = 0$. By a theorem of Serre and Cartan mentioned above (see [F-F3.25])  there exists such a number $s(2i)$ for all $i$. Let $p(N, k)$ be the smallest integer, such that for all $1\le i
\le k$ we have
$$ p(N,k) = \alpha (N, k, i) \cdot A_i\cdot s(2i)\cdot \beta (N, k, i)$$
for some positive integer $\beta (N, k, i)$. We shall construct a
map $ T: Gr _k (\C^N) \to Gr_k (\C^{2k(N-k)})$ such that
$$T^* (c_k ) = p(N,k) \cdot l_k \cdot c_k.\leqno(2.6)$$
Then,  taking into account of the  functoriality of the Chern classes,  the bundle $\hat E ^k$  defined by
$ \hat E ^k  = (f\circ  T )^* \gamma_k$  satisfies  the condition of Theorem 1.1.
Our map $T$ is  the composition of
$$ (f_{k,1} ^N, f_{k,2} ^N, \cdots , f_{k,k}^N) $$
where
$$f_{k,i}^N = G^N_{k,i} \circ \overline{F^{2i}_{l_i \cdot \beta(N,k,i)}} \circ \lambda_{k,i} ^N,$$
where  $  \overline{F^{2i}_{l_i \cdot s(2i)\cdot \beta(N,k,i)}}$ denotes the
restriction  of  the map $F ^{2i} _{l_i\cdot s(2i)\cdot \beta(N,k,i)}$ to
$K^{2k(N-k)} (\Z, 2i)$ (see Lemma 2.1). Because of our choice of $s(2i)$ and taking into account of Lemma 2.1, the map $T$ satisfies the condition (2.6).
 \QED
\medskip

\section {K\"ahler weak equivalence}

In this section we discuss some problems which arise in extending  the results in the previous section to the category of holomorphic bundles over complex or projective algebraic manifolds.

We would like to show another necessity for the notion of K\"ahler weak equivalence notion. 
Let $E^k$ be a complex  vector  bundle over a complex manifold $M^n$ and
$[f_{E^k}]$ be a the homotopy class of a classifying map $M \to Gr_k (\C^N)$ for $E^k$.
Clearly if $[f_{E^k}]$ contains a holomorphic map, then $E^k$ has a holomorphic
structure. But the converse statement is not true, because  the pull back  of
of any  positive $(1,1)$-cohomology classes via holomorphic map  is also  a non-positive
$(1,1)$-cohomology class. On the other hand there are many holomorphic vector bundles
whose first Chern class is a negative $(1,1)$-class. We shall say that a holomorphic vector bundle is {\bf positive}, if its classifying class contains a holomorphic map.

\medskip

{\bf 3.1. Lemma.} {\it Suppose that $M^m$ is a projective algebraic  manifold
and $E^k$ is a holomorphic vector bundle over $M^m$. Then $E^k$ is K\"ahler weakly equivalent to a positive holomorphic vector bundle.}

\medskip

{\it Proof.} This Lemma is a consequence of an well-known fact (see e.g. [G-H1978, Chapter 1, \S 5])
that a tensor of $E^k$ with a  some power $L^{\otimes l}$ of a K\"ahler line bundle $L$ admits enough holomorphic sections which serve as  a holomorphic map from $M^m$ to $Gr_k (\C^N)$,
where $\C^N$ is a subspace in $H^0 (M^m, \Oo (E^K))$. Furthermore, this holomorphic map is a classifying map for the holomorphic bundle $E\otimes L^{\otimes k}$, see
e.g. [G-H1978, Chapter 3, \S 3]. \QED

\medskip

We are lead by Lemma 3.1 to study the space $Hol(M^m, BU_k)$ of holomorphic maps from $M^m \to BU_k$. Denote by $Hodge (M)$ the group $H^{p,p} (M,\Q)$ and by
$[Hodge (M)]$ the quotient class $Hodge (M)/ \Q$ by multiplication.
We define then a map $C: Hol (M^m, BU_k)\to [Hodge ^k (M)]$ by $ C(f) = [f ^* (c_k)]$.
Then the Hodge conjecture is true, if and only if the image of map $C$ contains some neighborhood of a point $[P]\in [Hodge^k(M)]$ for some $P$ being  a power of a K\"ahler  class.
The problem in this naive   thinking  is that, $C$ maps a connected component  of $Hol(M^m,BU_k)$  to one point. It seems that we need to work every thing (including the Hodge theory) from the beginning in the  field of rationals.   Another possible way to do with  is mentioned in the introduction.

\medskip

\section{Appendix:  A proof of a theorem of Thom.}

Let $M ^n$ be an orientable differentiable manifold.

\medskip

{\bf A.1. Theorem} [Thom1954 Theorem II.25]. {\it For each
cohomology class $z \in H^k (M^{n}, \Z)$ there exists a number
$N(k, n)$ such that the class $N(k,n)\cdot z$ is the Euler class of an orientable vector
bundle on $M^n$. If $k= 2l$, then there exists a
number $N_1 (k,n)\ge N(k,n)$ such that the class $N_1 (k, n)\cdot
z$ is a top Chern class of a complex vector bundle on $M^n$. }

\medskip

Thom gave a detailed proof of Theorem A.1 for $G =  SO (k)$. He
noticed that his proof also works for $G = U(k)$ or $Sp (k)$.
Since we use Thom's theorem for $G = U(k)$ in 
[Le2005] as well as for our statement in the introduction  on the relation with Atiyah-Bott theorem, we feel a need for a detailed proof of Thom's theorem A.1
in this case.

\medskip

{\it Proof of Theorem A.1.} Suppose
that  $u \in H^{2k} (M^m, \Z)$. Then there is a map
$$f : M^m \to  K(\Z, 2k)$$
 such that $f^* (\tau^{2k}) = u$, where
$\tau^{2k}$ is the fundamental class of $H^k ( K (\Z, 2k), \Z)$.  Moreover we can assume that $f(M^m) \subset K^m (\Z, 2k)$, where $K^q(\Z, 2k)$ is the q-skeleton of the 
Eilenberg-McLane space $K(\Z, 2k)$.  
To prove Theorem A.1 it suffices to find a map
$$ h : K ^m (\Z, 2k) \to B_{U(k)}$$
such that for some positive number $N_1(k,m)$ we have
$$ h^* ( c_k) = N_1 (k,m) j ^* (\tau^{2k}),\leqno (2.2)$$
where $e_{2k} = c_k$ is the top Chern class of the universal bundle $V_{U(k)}$ over
$B_{U(k)}$ and $j$ is the embedding  $K^m (\Z, 2k) \to K (\Z , 2k)$. 

To find a map $h$ we apply Proposition 2.1. The main issue is to verify that
 the space $BU_k$ satisfies the condition for the space  $Y$ in Proposition 2.1. We use the same argument as that in our  proof of Lemma 2.3, actually the case
of $BU_k$ is easier, since the related exact sequences are simpler, The required
map $h$ can be constrcuted  in the same way as we did in our proof of Lemma 2.2 .)\QED

\medskip

{\bf Acknowledgement.}  This note has been written during my stay
at the Max-Planck-Institute f\"ur Mathematik in Leipzig. I am
grateful to J\"urgen Jost for  hospitality and  financial
support. I thank Dietmar Salamon for his interest in this project.

\bigskip

{\it email: hvle@math.cas.cz, \: hvle@mis.mpg.de}

\end{document}